\documentclass[11pt]{amsart}
\usepackage[top =2.2cm, bottom=2cm, right=3.5cm, left=3.5cm]{geometry}
\usepackage{amsmath,amsfonts,amsthm,amssymb,amscd}
\usepackage{tikz}
\usepackage{booktabs}
\usetikzlibrary{fadings}
\usetikzlibrary{patterns}
\usetikzlibrary{shadows.blur}
\usetikzlibrary{shapes}
\usepackage{float,graphicx}
\usepackage{newpxtext,newpxmath}
\usepackage{upgreek}
\usepackage{diagbox}
\usepackage{xcolor}
\usepackage{enumitem}
\usepackage[hidelinks]{hyperref}

\newcommand{\cN}{\mathcal{N}}

\newcommand{\cH}{\mathcal{H}}
\newcommand{\cT}{\mathcal{T}}
\newcommand{\cE}{\mathcal{E}}

\newcommand{\sfh}{\mathsf{h}}
\newcommand{\sfb}{\mathsf{b}}
\newcommand{\sfn}{\mathsf{n}}

\newcommand{\dist}{\operatorname{dist}}
\newtheorem{theorem}{Theorem}[section]

\theoremstyle{remark}
\newtheorem{remark}[theorem]{Remark}
\numberwithin{equation}{section}

\title{Exponentially Convergent Multiscale Finite Element Method}
\author{Yifan Chen, Thomas Y. Hou, Yixuan Wang}
\address{Applied and Computational Mathematics, Caltech, 91106}
\email{yifanc@caltech.edu, hou@cms.caltech.edu, roywang@caltech.edu}
\date{December 2022}
\subjclass[2010]{65N12, 65N15, 65N30, 31A35}
\begin{document}

\maketitle
\begin{center}
\textit{\small Dedicated to Professor Stanley Osher's 80th birthday with admiration and friendship}
\end{center}
\begin{abstract}
We provide a concise review of the exponentially convergent multiscale finite element method (ExpMsFEM) for efficient model reduction of PDEs in heterogeneous media without scale separation and in high-frequency wave propagation. ExpMsFEM is built on the non-overlapped domain decomposition in the classical MsFEM while enriching the approximation space systematically to achieve a nearly exponential convergence rate regarding the number of basis functions. Unlike most generalizations of MsFEM in the literature, ExpMsFEM does not rely on any partition of unity functions. 

In general, it is necessary to use function representations dependent on the right-hand side to break the algebraic Kolmogorov $n$-width barrier to achieve exponential convergence. Indeed, there are online and offline parts in the function representation provided by ExpMsFEM. The online part depends on the right-hand side locally and can be computed in parallel efficiently. 
The offline part contains basis functions that are used in the Galerkin method to assemble the stiffness matrix; they are all independent of the right-hand side, so the stiffness matrix can be used repeatedly in multi-query scenarios.

\end{abstract}
\section{Introduction}
Multiscale methods provide an efficient way to solve challenging PDEs. A few local basis functions adapted to the problem are constructed offline to provide an effective model reduction of the equation. One can then use the reduced model to compute the solution online, possibly with different right-hand sides and in a way much faster than solving the original equation. This property is beneficial in multi-query scenarios such as optimal design and inverse problems. Moreover, multiscale methods are inevitable for challenging problems in rough media and high-frequency wave propagation since standard numerical methods suffer from a vast number of degrees of freedom. See examples of the failure of finite element methods (FEMs) in elliptic equations with rough coefficients \cite{babuvska2000can} and the pollution effect in the Helmholtz equation \cite{babuska1997pollution}.

In this paper, we present the framework of ExpMsFEM, the exponentially convergent multiscale finite element method. It is a generalization of the classical MsFEM \cite{hou_multiscale_1997}. The main contribution of ExpMsFEM is the systematic improvement over MsFEM to achieve exponentially convergent accuracy regarding the number of basis functions. Also, unlike most generalizations of MsFEM in the literature, ExpMsFEM does not rely on the partition of unity functions to connect local and global approximation spaces. Instead, ExpMsFEM uses edge localization and coupling intrinsic to the non-overlapped domain decomposition to communicate the local and global approximations.

In the literature, exponentially convergent multiscale methods have been pioneered in the work of optimal basis \cite{babuska2011optimal} based on the partition of unity functions; see also the developments in \cite{smetana2016optimal,buhr2018randomized,chen2020randomized,babuvska2020multiscale, schleuss2020optimal, ma2021error, ma2021novel}. The work demonstrates the importance of Caccioppoli's inequality in establishing exponential convergence; more precisely, the inequality implies the \textit{low approximation complexity} of the restriction operator acting on harmonic-type functions. The theory of ExpMsFEM is also based on some arguments using Caccioppoli's inequality. Additionally, since no partition of unity functions is used, technical tools such as $C^{\alpha}$ estimates and trace theorems are needed to analyze ExpMsFEM. We will comment on the similarity and differences between the optimal basis work and ExpMsFEM at the end of the article.

This review is based on our previous work on exponentially convergent multiscale methods for elliptic equations \cite{chen2020exponential} and Helmholtz equations \cite{chen2021exponentially}. We focus on articulating the main ideas and the computational framework in the case of 2D stationary problems with homogeneous boundary data. We provide references for the detailed analysis in corresponding papers.

\subsection*{Organization} In Section \ref{sec: model problem}, we present the model problem that is the focus of this article. In Section \ref{sec: The ExpMsFEM Framework}, we present the motivation and framework of the ExpMsFEM. We provide numerical experiments to demonstrate the effectiveness of the ExpMsFEM framework in Section \ref{sec: Numerical Experiments}. In Section \ref{sec: Discussions}, we discuss related literature, future possibilities, and open questions. 

\section{Model Problem}
\label{sec: model problem}
Consider the model problem in a bounded domain $\Omega \subset \mathbb{R}^d$ with a Lipschitz boundary $\Gamma$. Here, $d=2$. For generality, the boundary can contain disjoint parts $\Gamma = \Gamma_1\cup \Gamma_2$ where $\Gamma_1$ corresponds to the Dirichlet boundary conditions and $\Gamma_2$ corresponds to the Neumann and Robin type boundary conditions. 
The model equation is: 
\begin{equation}
\label{eqn: model}
\left\{
\begin{aligned}
-\nabla \cdot(A\nabla u)+Vu&=f, \ \text{in} \ \Omega\\
u&=0, \ \text{on} \ \Gamma_1\\
A\nabla u\cdot\nu&=\beta u, \ \text{on} \  \Gamma_2 \, .
\end{aligned}
\right.
\end{equation}
Here, $A, V,\beta$ are functions in $L^{\infty}(\Omega)$ and can be rough, which makes the solution oscillating and difficult to solve. The vector $\nu$ is the outer normal to the boundary. 

In particular, when $V=0$, the equation is the standard elliptic equation \cite{chen2020exponential}. If $Vu=-k^2u$ and $u$ is a complex-valued function, one obtains the Helmholtz equation \cite{chen2021exponentially} with wavenumber $k$. 

The weak formulation of \eqref{eqn: model} is given by
\begin{equation}
\label{eqn: weak form}
a(u, v) :=(A\nabla u, \nabla {v})_{\Omega}+(V u,  {v})_{\Omega}  -( \beta u,  {v})_{\Gamma_2} =( f,  {v})_{\Omega} , \quad \forall v \in \cH(\Omega)\, ,
\end{equation}
where $(\cdot,\cdot)_X$ is the standard $L^2$ inner product on the set $X$.
The space for $v$ is $\cH(\Omega):=\{w \in H^1(\Omega): w|_{\Gamma_1}=0 \}$ and the solution $u \in \cH(\Omega)$. The energy norm $\|\cdot\|_{\cH(\Omega)}$ is defined as 
\[\|w\|_{\cH(\Omega)}^2:=(A\nabla w, \nabla {w})_{\Omega}+|(V w,  {w})_{\Omega}|\,.\]
Here, we adopt an abuse of notation that the space can be real-valued or complex-valued, depending on the context.

A generic assumption for $A$ is $0<A_{\min}\leq A(x)\leq A_{\max} < \infty$. We will present more detailed assumptions on $V,\beta$ later in specific problems that our theory in \cite{chen2020exponential,chen2021exponentially} covers. 
Indeed, the theory can encompass the case for very general $V$, provided that $|(Vu,u)|_{\Omega}\leq V_0(u,u)_{\Omega}$ for some constant $V_0$ and the PDE satisfies good stability estimates; see for example the rough Helmholtz example in \cite{chen2021exponentially}.
In this review, we mainly focus on the \textit{conceptual algorithmic framework} of solving the equation \eqref{eqn: model} via ExpMsFEM rather than a detailed analysis of the equation and the method.

\section{The ExpMsFEM Framework}
\label{sec: The ExpMsFEM Framework}
In subsection \ref{sec: Solving PDEs as function approximation}, we discuss the general recipe for solving PDEs as a function approximation problem. This motivates us to find accurate function representations to be used in the Galerkin method. We explain how ExpMsFEM manages to get exponentially convergent representations in subsections \ref{sec: Harmonic-bubble splitting}, \ref{sec: edge loc}, \ref{sec: Exponentially convergent SVD} and \ref{sec: alg ExpMsFEM}.
\subsection{Solving PDEs as function approximation}
\label{sec: Solving PDEs as function approximation}
By the standard finite element theory (e.g., \cite{brenner2008mathematical}), when using the Galerkin method to solve \eqref{eqn: weak form}, a key step is to find a function representation, or a space of basis functions that can approximate the solution accurately. More precisely, suppose the space is $S$, then, one usually wants
 \begin{equation}
    \label{apro proxy}
        {\eta(S)}:=\sup _{f \in L^{2}(\Omega) \backslash\{0\}} \inf _{v \in S} \frac{\left\|N(f)-v\right\|_{\mathcal{H}(\Omega)}}{\|f\|_{L^{2}(\Omega)}}
    \end{equation}
to be small. Here, $N: f \to u$ is the solution operator\footnote{Sometimes, $N$ is chosen to be the solution operator of the adjoint equation; for example see \cite{melenk2010convergence}.} of \eqref{eqn: model}.

For example, consider the elliptic equation with $V=0$ and $\Gamma_2=\emptyset$. In such case, the Galerkin method provides an optimal approximation of the solution in the space of basis functions with respect to the energy norm \cite{brenner2008mathematical, chen2020exponential}, due to the Galerkin orthogonality. Therefore, a small $\eta(S)$ directly implies a small error in the solution. For the Helmholtz equation, similar arguments hold based on the G\aa rding-type inequality, which leads to the quasi-optimality of the solution; see, for example, \cite{melenk2010convergence, chen2021exponentially}. The failure of many finite element methods in elliptic equations with rough coefficients \cite{babuvska2000can} and Helmholtz's equations \cite{babuska1997pollution} is due to the poor approximation property. $\eta(S)$ is typically not small if $S$ is the standard finite element space, such as the space of tent functions.

Conceptually, ExpMsFEM finds an exponentially convergent function representation of the solution through the following three steps: (1) harmonic-bubble splitting, (2) edge localization, (3) oversampling and exponentially convergent singular value decomposition (SVD). We will detail the three steps and discuss relevant rigorous results at the end of subsections \ref{sec: Harmonic-bubble splitting}, \ref{sec: edge loc}, and \ref{sec: Exponentially convergent SVD}. Then, we summarize the algorithm in subsection \ref{sec: alg ExpMsFEM}.
\subsection{Harmonic-bubble splitting}
\label{sec: Harmonic-bubble splitting}
Consider a shape regular and uniform partition of the domain $ \Omega $ into finite elements with a mesh size $H$. The collection of elements is denoted by $\cT_H=\{T_1, T_2,..., T_r\}$. Let $\cE_H=\{e_1,e_2,...,e_q\}$ be the collection of edges in the interior of $\Omega$. We use $\cN_H=\{x_1,x_2,...,x_p\}$ to denote the collection of interior nodes. We also use $E_H$ to denote the collection of interior edges as a set, i.e., $E_H=\bigcup_{e \in \cE_H} e \subset \Omega$. A more detailed explanation of the mesh structure can be found in \cite{chen2020exponential,chen2021exponentially}.

In each element $T\in\cT_H$,
we decompose the solution $u$ into $u=u_{T}^{\sfh}+u_T^{\sfb}$ such that
\begin{equation}
\label{eqn:decomposed}
\begin{aligned}
    &\left\{
    \begin{aligned}
    -\nabla \cdot (A \nabla u_T^\sfh )+V u^\sfh_T&=0, \ \text{in} \  T\\
    u_T^\sfh&=u, \ \text{on} \  \partial T \setminus (\Gamma_1\cup\Gamma_2)\\
    u_T^\sfh&=0, \ \text{on} \  \partial T \cap \Gamma_1\\
    A\nabla u_T^\sfh\cdot\nu&=\beta u_T^\sfh,\  \text{on} \ \partial T \cap \Gamma_2\, ,
    \end{aligned}
    \right.
    \\
    &\left\{
    \begin{aligned}
    -\nabla \cdot (A \nabla u^\sfb_T )+V u^\sfb_T&=f,\  \text{in} \  T\\
    u^\sfb_T&=0, \ \text{on} \  \partial T\setminus (\Gamma_1\cup\Gamma_2)\\
    u^\sfb_T&=0, \ \text{on} \  \partial T\cap \Gamma_1\\
   A\nabla u^\sfb_T\cdot\nu&=
  \beta u^\sfb_T, \ \text{on} \ \partial T \cap \Gamma_2 \, .
    \end{aligned}
    \right.
    \end{aligned}
    \end{equation}
In short, $u^{\sfh}_T$ incorporates the interior boundary value of $u$ on the element, while $u^{\sfb}_T$ contains information of the right-hand side. All equations in \eqref{eqn:decomposed} should be understood in the standard weak sense as in \eqref{eqn: weak form}.

We can further define a global decomposition $u=u^{\sfh}+u^{\sfb}$, such that for each $T$, it holds that $u^{\sfh}(x)=u^{\sfh}_T(x)$, $u^{\sfb}(x)=u^{\sfb}_T(x)$ when $x \in T$. 
 Here, the component $u^{\sfh}_T$ (resp. $u^{\sfh}$) is called the local (resp. global) \textit{harmonic part}, $u^{\sfb}_T$ (resp. $u^{\sfb}$) is the local (resp. global) \textit{bubble part}, of the solution $u$. Here, the harmonic part $u^\sfh$ is not necessarily a harmonic function due to the existence of $A$ and $V$, but it has a similar low complexity property that a harmonic function has, due to the iterative argument of Caccioppoli's inequality first proposed in \cite{babuska2011optimal}. We will discuss this low complexity property in subsection \ref{sec: Exponentially convergent SVD}.
 
 Now, in the representation $u=u^{\sfh}+u^{\sfb}$, the part $u^{\sfb}$ can be directly computed by solving local problems in parallel since the local boundary conditions are all known. We are left to deal with the part $u^\sfh$.
\begin{remark}
\label{rmk: harmonic bubble splitting}
    We discuss several theoretical concerns and possible generalizations below:
\begin{itemize}[leftmargin=*]
    \item A sufficient condition for the local components in \eqref{eqn:decomposed} to be well-defined is that the operator $u \to -\nabla \cdot (A \nabla u) + Vu$ (as well as the corresponding boundary conditions) is elliptic in each local element, implied by the Poincar\'e inequality. In \cite{chen2020exponential}, we consider elliptic equations with $V = 0$ and $\Gamma_2 = \emptyset$, so this condition is satisfied. In \cite{chen2021exponentially}, we consider the Helmholtz equation where $V <0, |V|=O(k^2)$ and $\operatorname{Re}\beta =0, \operatorname{Im}\beta = O(k)$. For such a case, the elliptic property is guaranteed when $H = O(1/k)$.
    \item For the global components $u^\sfh$, $u^{\sfb}$ to be well-defined, we need the condition that the solution $u$ is continuous. This can be guaranteed by the $C^{\alpha}$ estimates of the equation \eqref{eqn: model} under the assumptions mentioned earlier; see discussions in \cite{chen2020exponential,chen2021exponentially}.
    \item We can generalize the above decomposition to PDEs with inhomogeneous boundary conditions. To achieve so, we incorporate these boundary data into the equation for $u^\sfb$; see also Section 5.3 in \cite{chen2021exponentially} for a concrete example of problem with inhomogeneous boundary data.
\end{itemize}
\end{remark}

\subsection{Edge localization} 
\label{sec: edge loc}
The next step is to find some local basis functions that accurately approximate $u^\sfh$. ExpMsFEM uses the idea of \textit{edge localization} to localize this approximation task.

First, we define the ``harmonic extension'' operator $Q_{E_H}$ that maps the edge values $\tilde{u}^\sfh = u^\sfh|_{E_H} \in H^{1/2}(E_H)$ to $ u^\sfh \in H^1(\Omega)$, through the relation in the first set of equation in \eqref{eqn:decomposed}. Here, we adopt the convention that if we write a tilde on the top of a function, it is the restriction of this function on the edge set. We have that $u^\sfh = Q_{E_H}\tilde{u}^\sfh = Q_{E_H}\tilde{u}$, since $u^\sfh$ and $u$ have the same edge values.

Then, let $C(E_H)$ be the space of continuous functions on $E_H$. We consider the edge interpolation operator $I_H: H^{1/2}(E_H) \cap C(E_H) \to H^{1/2}(E_H) \cap C(E_H)$ such that \[I_H \tilde{u} = \sum_{x_i \in \cN_H} \tilde{u}(x_i)\tilde{\psi}_i\] where the edge function $\tilde{\psi}_i$ is linear on $E_H$ and satisfies $\tilde{\psi}_i(x_j) = \delta_{ij}$. Note that by the convention of our notation we have $\psi_i = Q_{E_H}\tilde{\psi}_i \in H^1(\Omega)$. It is worth noting that $\psi_i's$ are the basis functions used in the vanilla MsFEM.

With the interpolation operator, we can write 
\[Q_{E_H}\tilde{u} = Q_{E_H}(\tilde{u}-I_H\tilde{u}) + \sum_{x_i \in \cN_H} u(x_i) \psi_i\, . \]

Now, the residue $\tilde{u}-I_H\tilde{u}$ is zero at each interior node. This property allows us to localize the residue to each edge. Indeed, by an abuse of notation, we can write
\begin{equation}
    Q_{E_H}(\tilde{u}-I_H\tilde{u}) = \sum_{e \in \cE_H} Q_{E_H}(\tilde{u}-I_H\tilde{u})|_{e}\, ,
\end{equation}
where we equate the function $(\tilde{u}-I_H\tilde{u})|_{e}$ that is defined on $e$ to its zero extension to $E_H$, so that $(\tilde{u}-I_H\tilde{u})|_{e}\in H^{1/2}(E_H)$ and thus $Q_{E_H}(\tilde{u}-I_H\tilde{u})|_{e}$ makes sense.

Therefore, we localize the approximation task of $u^\sfh$ to $Q_{E_H}(\tilde{u}-I_H\tilde{u})|_{e}$, which is defined for each edge $e$.

\begin{remark}
\label{rmk: theory for edge localization}
    Again, we discuss several theoretical concerns below:
    \begin{itemize}[leftmargin=*]
        \item Once the condition in Remark \ref{rmk: harmonic bubble splitting} is satisfied, the extension operator $Q_{E_H}$ is well-defined because the local equation is elliptic.
        \item According to the comment in Remark \ref{rmk: harmonic bubble splitting}, the solution $u$ is continuous, so the nodal interpolation $I_H \tilde{u}$ is well-defined.
        \item One can rigorously show that if we can approximate each local term with \[\|Q_{E_H}(\tilde{u}-I_H\tilde{u})|_{e} - w_e\|_{\cH(\Omega)} \leq \epsilon_e \, ,\]
        then the global approximation error satisfies
        \[\|Q_{E_H}(\tilde{u}-I_H\tilde{u}) - \sum_{e \in \cE_H} w_e\|_{\cH(\Omega)}^2 \leq C_{\mathrm{mesh}}\sum_{e \in \cE_H}\epsilon_e^2\, , \]
        where $C_{\mathrm{mesh}}$ is a constant dependent on the mesh structure only. In our previous work \cite{chen2020exponential,chen2021exponentially}, we formalize the approximation in the edge space via the $H_{00}^{1/2}(e)$ norm, which is equivalent to the $\cH(\Omega)$ norm here after the extension by $Q_{E_H}$; see Proposition 2.5 and Theorem 2.6 in \cite{chen2020exponential}. In this review paper, we explain the ideas using $Q_{E_H}$ rather than $H_{00}^{1/2}(e)$, since the former is more concise in an algorithm-focused exposition.
        
        We call the step from local approximation to global approximation \textit{edge coupling}.
    \end{itemize}
\end{remark}
\subsection{Exponentially convergent SVD}
\label{sec: Exponentially convergent SVD}
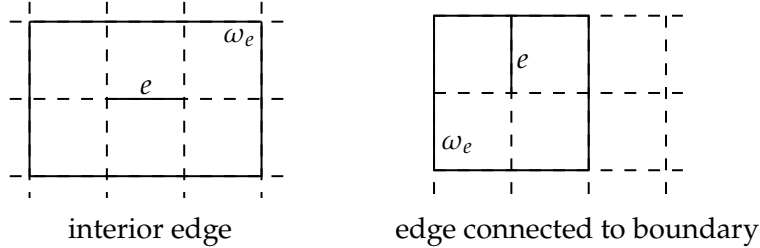
\begin{figure}[ht]
        \centering
\tikzset{every picture/.style={line width=0.75pt}} 
\begin{tikzpicture}[x=0.75pt,y=0.75pt,yscale=-1,xscale=1]

\draw  [draw opacity=0][dash pattern={on 4.5pt off 4.5pt}] (60,64) -- (201.5,64) -- (201.5,163) -- (60,163) -- cycle ; \draw  [dash pattern={on 4.5pt off 4.5pt}] (70,64) -- (70,163)(109,64) -- (109,163)(148,64) -- (148,163)(187,64) -- (187,163) ; \draw  [dash pattern={on 4.5pt off 4.5pt}] (60,74) -- (201.5,74)(60,113) -- (201.5,113)(60,152) -- (201.5,152) ; \draw  [dash pattern={on 4.5pt off 4.5pt}]  ;
\draw    (109,113) -- (148,113) ;
\draw    (70,74) -- (187,74) ;
\draw    (187,74) -- (187,152) -- (70,152) -- (70,74) ;
\draw  [draw opacity=0][dash pattern={on 4.5pt off 4.5pt}] (274,71) -- (399.5,71) -- (399.5,162) -- (274,162) -- cycle ; \draw  [dash pattern={on 4.5pt off 4.5pt}] (274,71) -- (274,162)(313,71) -- (313,162)(352,71) -- (352,162)(391,71) -- (391,162) ; \draw  [dash pattern={on 4.5pt off 4.5pt}] (274,71) -- (399.5,71)(274,110) -- (399.5,110)(274,149) -- (399.5,149) ; \draw  [dash pattern={on 4.5pt off 4.5pt}]  ;
\draw    (313,71) -- (313,110) ;
\draw    (274,71) -- (274,149) -- (352,149) -- (352,71) -- cycle ;

\draw (124,102) node [anchor=north west][inner sep=0.75pt]   [align=left] {$\displaystyle e$};
\draw (166,76) node [anchor=north west][inner sep=0.75pt]   [align=left] {$\displaystyle \omega _{e}$};
\draw (88,171) node [anchor=north west][inner sep=0.75pt]   [align=left] {interior edge};
\draw (253,171) node [anchor=north west][inner sep=0.75pt]   [align=left] {edge connected to boundary};
\draw (314,89) node [anchor=north west][inner sep=0.75pt]   [align=left] {$\displaystyle e$};
\draw (276,130) node [anchor=north west][inner sep=0.75pt]   [align=left] {$\displaystyle \omega _{e}$};
\end{tikzpicture}
        \caption{Illustration of oversampling domains. On the right, we use an edge connected to the upper boundary as an illustrating example.}
        \label{fig:os domain}
    \end{figure}

Recall that by using the harmonic-bubble splitting and edge localization, we get the representation
\begin{equation}
    u = u^{\sfh}+u^{\sfb} = \sum_{e \in \cE_H} Q_{E_H}(\tilde{u}-I_H\tilde{u})|_{e} + \sum_{x_i \in \cN_H} u(x_i) \psi_i + u^{\sfb}\, .
\end{equation}
ExpMsFEM then relies on oversampling and local SVD to get an exponentially convergent approximation of each $Q_{E_H}(\tilde{u}-I_H\tilde{u})|_{e}$. For each $e$, consider an oversampling domain $w_e \supset e$. Any domain containing $e$ in the interior may be used, and as an illustrative example, we set 
\begin{equation*}
    \label{eqn: os domain 1 layer}
        \omega_e=\overline{\bigcup \{T\in \cT_H: \overline{T} \cap e \neq \emptyset\}}\, .
    \end{equation*} 
An illustration of this choice for a quadrilateral mesh is given in Figure \ref{fig:os domain}. 

We can view $(\tilde{u}-I_H\tilde{u})|_{e}$ as the image of an operator acting on $u|_{\omega_e} \in H^1(\omega_e)$. We denote this operator by $R_e$ such that $Q_{E_H}(\tilde{u}-I_H\tilde{u})|_{e} = Q_{E_H}R_e (u|_{\omega_e})$. Now, we apply the harmonic-bubble splitting in subsection \ref{sec: Harmonic-bubble splitting} to the domain $\omega_e$, which leads to $u|_{\omega_e} = u_{\omega_e}^\sfh + u_{\omega_e}^\sfb$. It follows that
\begin{equation}
    Q_{E_H}(\tilde{u}-I_H\tilde{u})|_{e} = Q_{E_H}R_e u^\sfh_{\omega_e} + Q_{E_H}R_e u^\sfb_{\omega_e}\, .
\end{equation}
The term $R_e u^\sfh_{\omega_e}$ is a restriction of a harmonic part. As we mentioned at the beginning of this article, one can prove that the restriction operator acting on harmonic-type functions is of \textit{low approximation complexity}. More precisely, consider the space of harmonic parts in $\omega_e$, defined via 
\begin{equation}
    \begin{aligned}
        U(\omega_e):= \{v \in \cH(\omega_e): &-\nabla \cdot (A \nabla v)+Vv=0, \text{ in } \omega_e\\
        &A\nabla v\cdot\nu=\beta v, \text{ on } \Gamma_1 \cap \partial \omega_e\} 
        \, .
    \end{aligned}
    \end{equation}
The space is equipped with the norm $\|\cdot\|_{\cH(\omega_e)}$. Then, one can show that the left singular values (in descending order) of the local operator 
\[Q_{E_H}R_e: (U(\omega_e),\|\cdot\|_{\cH(\omega_e)}) \to (\cH(\Omega), \|\cdot\|_{\cH(\Omega)})\]
decays as $\lambda_{e,m} \leq C\exp(-bm^{\frac{1}{d+1}})$ in dimension $d$, for some generic constant $C,b$ independent of $m$ and $H$. Equivalently, if we write the left singular vectors as $v_{e,m} \in H^1(\Omega)$, which is local and supported in the neighboring elements of the edge $e$, then there exists some coefficient $b_{e,j}$ such that
\begin{equation}
    \|Q_{E_H}R_e u^\sfh_{\omega_e} - \sum_{1\leq j \leq m} b_{e,j} v_{e,j} \|_{\cH(\Omega)} \leq C\exp(-bm^{\frac{1}{d+1}})\|u^\sfh_{\omega_e}\|_{\cH(\omega_e)}\, .
\end{equation}
For more details, see Theorem 3.10 in \cite{chen2021exponentially}. Then, summing these local errors up, we get 
\begin{equation}
\begin{aligned}
    \sum_{e \in \cE_H} \|u^\sfh_{\omega_e}\|_{\cH(\omega_e)}^2 &\leq 2\sum_{e \in \cE_H} (\|u|_{\omega_e}\|_{\cH(\omega_e)}^2 + \|u_{\omega_e}^\sfb\|_{\cH(\omega_e)}^2)\\
    & = O(\|u\|^2_{\cH(\Omega)}+\|f\|^2_{L^2(\Omega)})\, ,
\end{aligned}
\end{equation}
where we used the fact that $\|u_{\omega_e}^\sfb\|_{\cH(\omega_e)} = O(\|f\|_{L^2(\omega_e)})$ by the elliptic estimate.

Combining the above  estimates with edge coupling in Remark \ref{rmk: theory for edge localization}, we get the representation
\begin{equation}
\label{eqn: exp representation}
\begin{aligned}
    u = u^{\sfh}+u^{\sfb} = &\sum_{e \in \cE_H} \sum_{1\leq j \leq m} b_{e,j} v_{e,j} + \sum_{x_i \in \cN_H} u(x_i) \psi_i + u^{\sfn} \\
    &+ O\left(\exp(-bm^{\frac{1}{d+1}})(\|u\|_{\cH(\Omega)}+\|f\|_{L^2(\Omega)})\right)\, ,
\end{aligned}
\end{equation}
where $u^\sfn := u^\sfb + \sum_{e \in \cE_H} Q_{E_H}R_e u^\sfb_{\omega_e}$ is a part that depends on $f$ locally.
\begin{remark}
    We discuss several theoretical aspects and the implication of the above representation.
    \begin{itemize}[leftmargin=*]
        
        \item The proof of the exponentially decaying singular values of $Q_{E_H}R_e$ is based on two steps. The first step is the iterative argument of Caccioppoli's inequality, first proposed in \cite{babuska2011optimal} and then refined in \cite{ma2021novel}. It shows that the singular values of the restriction operator on $U(\omega_e)$, which restricts a function from the original domain $\omega_e$ to a subdomain $\omega^* \supset e$, decay nearly exponentially fast. The second step is based on a stability estimate of the operator $Q_{E_H}R_e$ acting on $U(\omega^*)$; see Lemma 3.10 in \cite{chen2020exponential} or Lemma 6.1, 6.2 in \cite{chen2021exponentially}. 
        \item We can understand that the oversampling technique is used to take advantage of the low complexity property of the restriction operator. Historically, the idea of oversampling was proposed in \cite{hou_multiscale_1997} to reduce the resonance error in MsFEM.
        \item The remarkable thing about the representation in \eqref{eqn: exp representation} is the exponentially decaying error bound. 
        
        First, for elliptic equations with rough coefficients, the error bound implies that these basis functions can capture the behavior of the solution, which is a hard task for FEMs. Therefore, ExpMsFEM overcomes the difficulty of rough coefficients.
        
        Second, for the Helmholtz equation, the stability constant of the solution operator can depend on $k$; indeed, this is the main cause of the pollution effect \cite{babuska1997pollution}. Denote the stability constant by $C_{\text{stab}}(k)$ such that $\|u\|_{\cH(\Omega)} \leq C_{\text{stab}}(k)\|f\|_{L^2(\Omega)}$. A prevalent and reasonable assumption on the constant is that of polynomial growth, namely $C_{\mathrm{stab}}(k)\leq C(1+k^\gamma)$ for some constants $\gamma$ and $C$; see, for example, \cite{lafontaine2019most}. In such case, we can further bound the error by
        \[\exp(-bm^{\frac{1}{d+1}})(\|u\|_{\cH(\Omega)}+\|f\|_{L^2(\Omega)}) \leq  \exp(-bm^{\frac{1}{d+1}})(C(1+k^\gamma)+1)\|f\|_{L^2(\Omega)}\, .\]
        Therefore, once the number of basis functions per edge $m \sim \log^{d+1}(k)$ (logarithmically on $k$ only), the approximation error can be uniformly small for all $k$. It implies that the quantity $\eta(S)$ in \eqref{apro proxy} is small, which is important in determining the error of Galerkin's methods. In this sense, ExpMsFEM overcomes the difficulty of the pollution effect by using basis functions whose number scales at most $\log^{d+1}(k)$.
        \item The exponentially accurate representation in \eqref{eqn: exp representation} will not be possible if we do not use terms dependent on the right-hand side. Indeed, using basis functions independent of $f$, the optimal approximation error rate will be algebraic if the right-hand side is in $L^2(\Omega)$ only, due to well-known results in approximation theory (the Kolmogorov $n$-width \cite{pinkus2012n, melenk2000n}); see also the complexity analysis of the Green function of Helmholtz's equation \cite{engquist2018approximate}. From this perspective, we can understand that ExpMsFEM breaks the Kolmogorov barrier by using \textit{nonlinear model reduction} \cite{peherstorfer2022breaking}, i.e., the basis functions can depend on the input of the model, here the right-hand side.
    \end{itemize}
\end{remark}
\subsection{The solver based on ExpMsFEM}
\label{sec: alg ExpMsFEM}
Now, we can use the representation in \eqref{eqn: exp representation} to solve the equation efficiently. First, we form $\psi_i, v_{e,j}$ by computing the local extension $Q_{E_H}\tilde{\psi}_i$ for each node and the top-$m$ left singular vectors $v_{e,j}, 1\leq j \leq m$ of the local operator $Q_{E_H}R_e$ for each $e$; problems on different nodes and edges are independent and parallelizable. These become our offline basis functions.

For any right-hand side $f$, we compute the online part $u^\sfn$ by solving local linear equations involving $f$. This step can be parallelized. 

Then, we form an effective equation for $u-u^\sfn$ as
\begin{equation}
\label{eqn: effective eqn for u h}
    a(u-u^\sfn,v)=(f,v)_{\Omega}-a(u^\sfn,v)\, ,
\end{equation}
for any $v \in \cH(\Omega)$. We solve the equation for $u-u^\sfn$ using a Galerkin method. As an example, using the Ritz-Galerkin method, we choose 
\[S = \mathrm{span}~\{\psi_i\ \text{for } x_i \in \cN_H, \ v_{e,j} \  \text{for } 1\leq j \leq m, e \in \cE_H \}\, ,\] 
and find a numerical solution $u_S \in S$
that satisfies
\begin{equation}
\label{eqn: Galerkin}
    a(u_S,v)=(f,v)_{\Omega}-a(u^\sfn,v)\, ,
\end{equation}
for any $v \in S$. The final numerical solution is given by $u_S + u^\sfn$. We call $u^\sfn$ the online part and $u_S$ the offline part since $u_S$ lies in a space that is independent of $f$.

Note that in the Galerkin method for solving $u_S$, the stiffness matrix only needs to be assembled once and can be used for different $f$ afterward. We can understand \eqref{eqn: effective eqn for u h} as a reduced model of the original equation.
\begin{remark}
    We discuss several theoretical aspects regarding the effectiveness of the above method.
    \begin{itemize}[leftmargin=*]
        \item The accuracy of the numerical solution is due to the quasi-optimality property mentioned earlier in subsection \ref{sec: Solving PDEs as function approximation}: once $\eta(S)$ is small, the solution error is of the same order compared to the optimal approximation using the basis functions, which is exponentially small according to the representation \eqref{eqn: exp representation}.
        \item When the solution is complex-valued, such as in the Helmholtz equations, we can use both the Ritz and Petrov versions of the Galerkin methods; for the former, if $\overline{S} \neq S$, we need to replace $S$ by $S+\overline{S}$; see discussions in \cite{chen2021exponentially}.
        \item One thing worth noting is that $\|u^\sfn\|_{\cH(\Omega)}$ is of order $O(H)$, due to the standard elliptic estimate \cite{chen2020exponential,chen2021exponentially}. Therefore, if we aim for $O(H)$ accuracy only, we can ignore this part, and simply setting $u^\sfn$ = 0 in the above algorithm will lead to a solution accurate up to $O(H)$.
    \end{itemize}
\end{remark}

\section{Numerical Experiments}
\label{sec: Numerical Experiments}
\begin{figure}[htbp]
\centering

\tikzset{every picture/.style={line width=0.75pt}} 

\begin{tikzpicture}[x=0.75pt,y=0.75pt,yscale=-1,xscale=1]
\draw  [draw opacity=0][dash pattern={on 4.5pt off 4.5pt}] (170.5,75) -- (365.5,75) -- (365.5,231) -- (170.5,231) -- cycle ; \draw  [color={rgb, 255:red, 0; green, 0; blue, 0 }  ,draw opacity=1 ][dash pattern={on 4.5pt off 4.5pt}] (209.5,75) -- (209.5,231)(248.5,75) -- (248.5,231)(287.5,75) -- (287.5,231)(326.5,75) -- (326.5,231) ; \draw  [color={rgb, 255:red, 0; green, 0; blue, 0 }  ,draw opacity=1 ][dash pattern={on 4.5pt off 4.5pt}] (170.5,114) -- (365.5,114)(170.5,153) -- (365.5,153)(170.5,192) -- (365.5,192) ; \draw  [color={rgb, 255:red, 0; green, 0; blue, 0 }  ,draw opacity=1 ][dash pattern={on 4.5pt off 4.5pt}] (170.5,75) -- (365.5,75) -- (365.5,231) -- (170.5,231) -- cycle ;
\draw  [draw opacity=0][dash pattern={on 0.84pt off 2.51pt}] (170.5,75) -- (365.5,75) -- (365.5,231) -- (170.5,231) -- cycle ; \draw  [color={rgb, 255:red, 0; green, 0; blue, 0 }  ,draw opacity=1 ][dash pattern={on 0.84pt off 2.51pt}] (183.5,75) -- (183.5,231)(196.5,75) -- (196.5,231)(209.5,75) -- (209.5,231)(222.5,75) -- (222.5,231)(235.5,75) -- (235.5,231)(248.5,75) -- (248.5,231)(261.5,75) -- (261.5,231)(274.5,75) -- (274.5,231)(287.5,75) -- (287.5,231)(300.5,75) -- (300.5,231)(313.5,75) -- (313.5,231)(326.5,75) -- (326.5,231)(339.5,75) -- (339.5,231)(352.5,75) -- (352.5,231) ; \draw  [color={rgb, 255:red, 0; green, 0; blue, 0 }  ,draw opacity=1 ][dash pattern={on 0.84pt off 2.51pt}] (170.5,88) -- (365.5,88)(170.5,101) -- (365.5,101)(170.5,114) -- (365.5,114)(170.5,127) -- (365.5,127)(170.5,140) -- (365.5,140)(170.5,153) --
(365.5,153)(170.5,166) -- (365.5,166)(170.5,179) -- (365.5,179)(170.5,192) -- (365.5,192)(170.5,205) -- (365.5,205)(170.5,218) -- (365.5,218) ; \draw  [color={rgb, 255:red, 0; green, 0; blue, 0 }  ,draw opacity=1 ][dash pattern={on 0.84pt off 2.51pt}] (170.5,75) -- (365.5,75) -- (365.5,231) -- (170.5,231) -- cycle ;
\draw    (209.5,114) -- (209.5,192) -- (326.5,192) -- (326.5,114) -- cycle ;

\draw [color={rgb, 255:red, 0; green, 0; blue, 0 }  ,draw opacity=1 ][line width=0.75]    (248.5,153) -- (287.5,153) ;
\draw  [draw opacity=0][dash pattern={on 0.84pt off 2.51pt}] (407.5,157) -- (507.5,157) -- (507.5,218) -- (407.5,218) -- cycle ; \draw  [color={rgb, 255:red, 0; green, 0; blue, 0 }  ,draw opacity=1 ][dash pattern={on 0.84pt off 2.51pt}] (417.5,157) -- (417.5,218)(430.5,157) -- (430.5,218)(443.5,157) -- (443.5,218)(456.5,157) -- (456.5,218)(469.5,157) -- (469.5,218)(482.5,157) -- (482.5,218)(495.5,157) -- (495.5,218) ; \draw  [color={rgb, 255:red, 0; green, 0; blue, 0 }  ,draw opacity=1 ][dash pattern={on 0.84pt off 2.51pt}] (407.5,167) -- (507.5,167)(407.5,180) -- (507.5,180)(407.5,193) -- (507.5,193)(407.5,206) -- (507.5,206) ; \draw  [color={rgb, 255:red, 0; green, 0; blue, 0 }  ,draw opacity=1 ][dash pattern={on 0.84pt off 2.51pt}]  ;
\draw  [draw opacity=0][dash pattern={on 4.5pt off 4.5pt}] (406.5,73) -- (508.5,73) -- (508.5,135) -- (406.5,135) -- cycle ; \draw  [color={rgb, 255:red, 0; green, 0; blue, 0 }  ,draw opacity=1 ][dash pattern={on 4.5pt off 4.5pt}] (416.5,73) -- (416.5,135)(455.5,73) -- (455.5,135)(494.5,73) -- (494.5,135) ; \draw  [color={rgb, 255:red, 0; green, 0; blue, 0 }  ,draw opacity=1 ][dash pattern={on 4.5pt off 4.5pt}] (406.5,83) -- (508.5,83)(406.5,122) -- (508.5,122) ; \draw  [color={rgb, 255:red, 0; green, 0; blue, 0 }  ,draw opacity=1 ][dash pattern={on 4.5pt off 4.5pt}]  ;

\draw (263.5,140) node [anchor=north west][inner sep=0.75pt]   [align=left] {$\displaystyle e$};
\draw (303.5,115) node [anchor=north west][inner sep=0.75pt]   [align=left] {$\displaystyle \omega _{e}$};
\draw (414,136) node [anchor=north west][inner sep=0.75pt]   [align=left] {coarse mesh};
\draw (424,218) node [anchor=north west][inner sep=0.75pt]   [align=left] {fine mesh};

\end{tikzpicture}

\caption{Two level mesh: a fraction}
\label{fig:mesh1}
\end{figure}
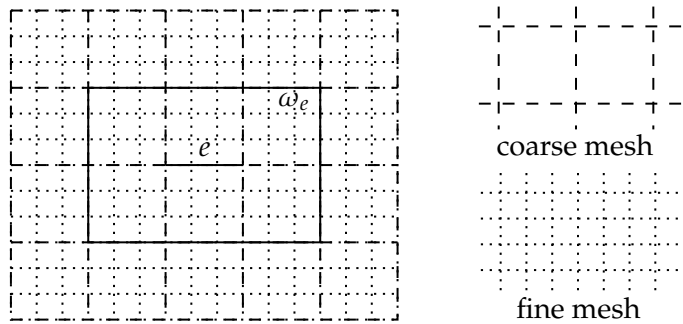

In this section, we present some numerical experiments to demonstrate the effectiveness of ExpMsFEM. For all the experiments, we consider the domain $\Omega=[0,1]\times [0,1]$ and discretize it by a uniform two-level quadrilateral mesh; see a fraction of this mesh in Figure \ref{fig:mesh1}, where we also show an edge $e$ and its oversampling domain $\omega_e$ in solid lines. 
 The coarse and fine mesh sizes are denoted by $H$ and $h$, respectively. 
 
For a given equation,  we compute the reference solution $u_{\text{ref}}$ using the classical FEM on the fine mesh with a sufficiently small $h$, which we choose to be $h=1/1024$. By \textit{a posteriori} estimates, we can check that the fine mesh indeed resolves the corresponding problems; thus, the associated fine mesh solutions could serve as accurate reference solutions for all of our numerical examples. In our numerical computation, we solve local problems that are required in the ExpMsFEM framework using the fine mesh. For detailed implementation, we refer to \cite{chen2020exponential,chen2021exponentially}.
 \begin{remark}[Accuracy on the discrete level]
For simplicity of presentation, we do not provide error analysis of ExpMsFEM on the fully discrete level, where the accuracy of the local problems can depend on the resolution of the fine grid. For a detailed error estimate on the fully discrete level in the context of partition of unity methods, see, for example, \cite{ma2021error, ma2021wavenumber}.
\end{remark}
The accuracy of a numerical solution $u_{\mathrm{sol}}$ is computed by comparing it with the reference solution $u_{\text{ref}}$ on the fine mesh. The accuracy will be measured both in the $L^2$ norm and energy norm: 
    \begin{equation}
\label{rel_error}
\begin{aligned}
e_{L^{2}}&=\frac{\|u_{\text{ref}}-u_{\mathrm{sol}}\|_{L^{2}(\Omega)}}{\|u_{\text{ref}}\|_{L^{2}(\Omega)}}\, ,\\
e_{\cH}&=\frac{\|u_{\text{ref}}-u_{\mathrm{sol}}\|_{\cH(\Omega)}}{\|u_{\text{ref}}\|_{\cH(\Omega)}}\, . 
\end{aligned}
\end{equation}

In subsection \ref{sec: exp period multiscale}, we consider an elliptic equation where the coefficient $A(x)$ is periodic but contains multiple scales. This example demonstrates the exponential accuracy of ExpMsFEM. In subsection \ref{sec: exp high contrast}, we consider an elliptic equation where $A(x)$ is of high contrast. This example shows the robustness of ExpMsFEM regarding the high contrast. In subsection \ref{sec: exp Helmholtz}, an instance of Helmholtz's equation with rough media and mixed boundary conditions is presented. This example illustrates the effectiveness of ExpMsFEM in solving general indefinite Helmholtz's equations.
\subsection{A periodic example with multiple spatial scales}
\label{sec: exp period multiscale}
In the first example, we consider an elliptic problem ($V=0$) with multiple spatial scales. We choose coefficient $A$ with five scales as follows: 
   \begin{equation}
\begin{aligned} 
A(x)=\frac{1}{6}\left(\frac{1.1+\sin \left(2 \pi x_1 / \epsilon_{1}\right)}{1.1+\sin \left(2 \pi x_2 / \epsilon_{1}\right)}+\frac{1.1+\sin \left(2 \pi x_2 / \epsilon_{2}\right)}{1.1+\cos \left(2 \pi x_1 / \epsilon_{2}\right)}+\frac{1.1+\cos \left(2 \pi x_1 / \epsilon_{3}\right)}{1.1+\sin \left(2 \pi x_2 / \epsilon_{3}\right)}\right.\\\left.+\frac{1.1+\sin \left(2 \pi x_2 / \epsilon_{4}\right)}{1.1+\cos \left(2 \pi x_1 / \epsilon_{4}\right)}+\frac{1.1+\cos \left(2 \pi x_1 / \epsilon_{5}\right)}{1.1+\sin \left(2 \pi x_2 / \epsilon_{5}\right)}+\sin \left(4 x_1^{2} x_2^{2}\right)+1\right) \, ,\end{aligned}
\end{equation}
where $x=(x_1,x_2)$, $\epsilon_1=1/5$, $\epsilon_2=1/13$, $\epsilon_3=1/17$, $\epsilon_4=1/31$, $\epsilon_5=1/65$. We choose homogeneous Dirichlet boundary conditions, i.e., $\Gamma_2=\emptyset$. We set $f=-1$.

\begin{figure}[ht]
    \centering
    \includegraphics[width=6cm]{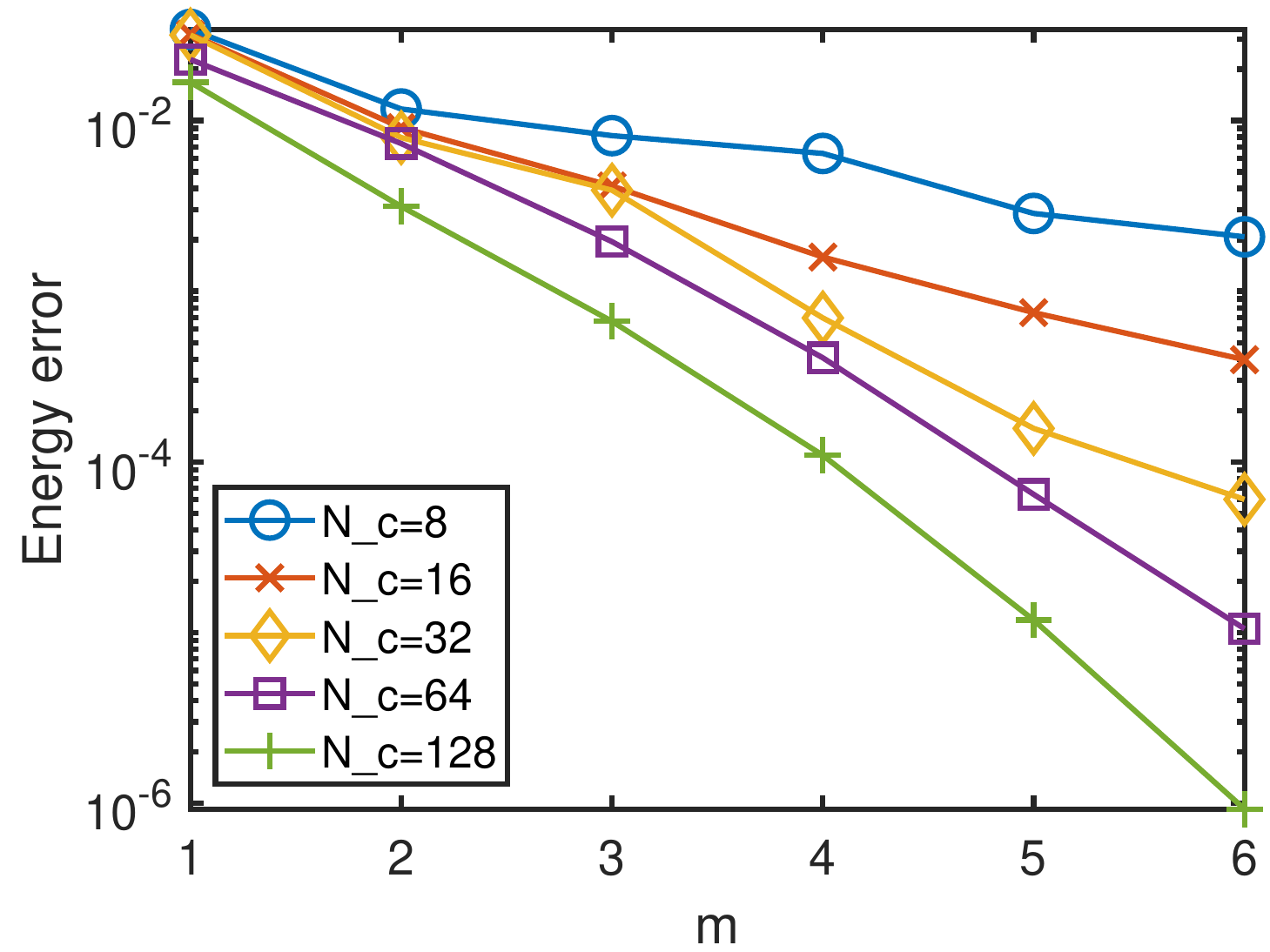}
    \includegraphics[width=6cm]{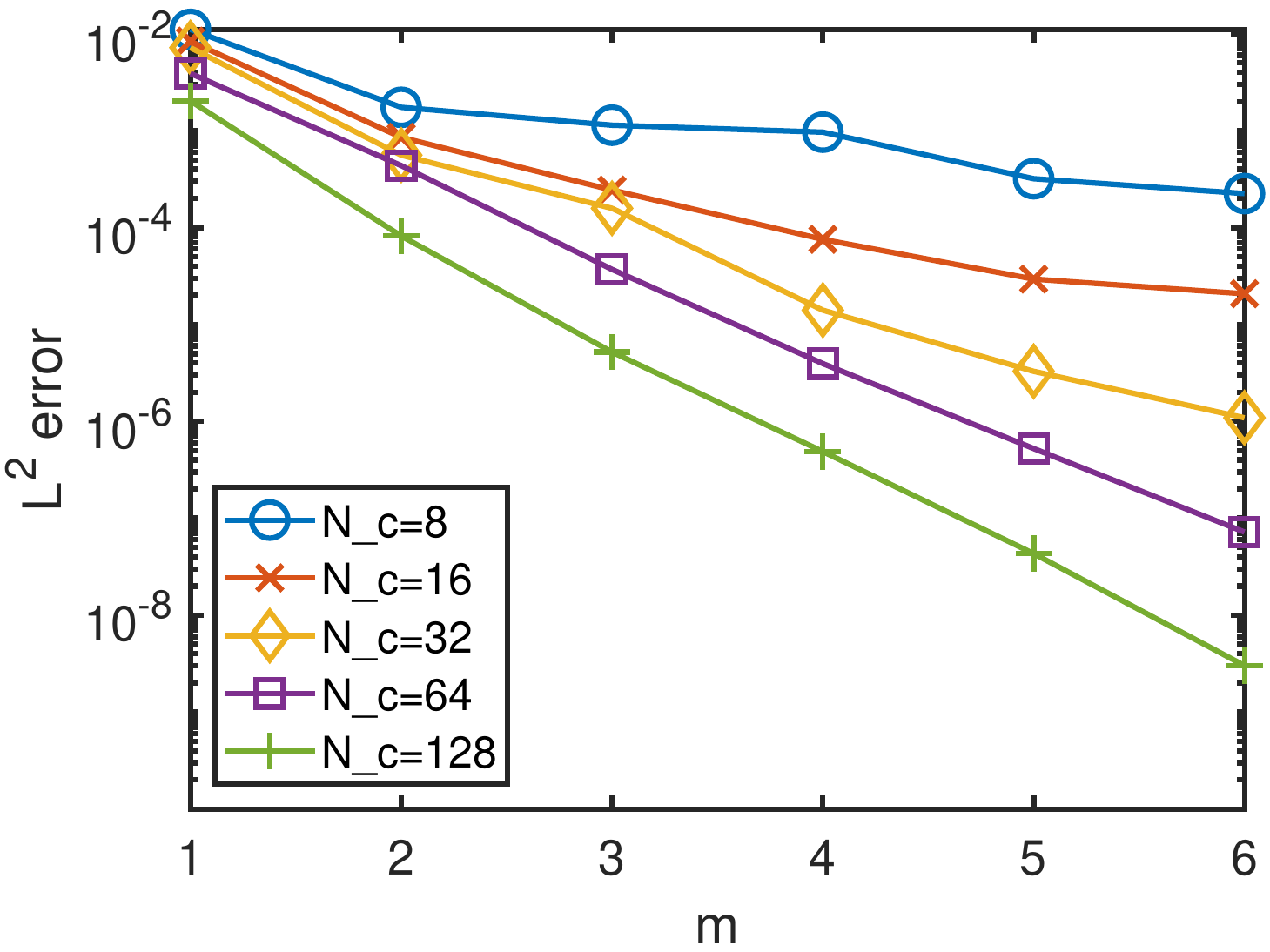}
    \caption{Numerical results for the periodic example. Left: $e_{\cH}$ versus $m$; right: $e_{L^2}$ versus $m$.}
    \label{fig:eg1}
\end{figure}

In this example, we illustrate the exponential accuracy and the convergence rate with respect to the coarse mesh size $H$. We take $H=2^{-i}$, $i=3,4,...,7$ and take $m=1,2,...,6$ for each $H$. The numerical results are shown in Figure \ref{fig:eg1}, where $N_c=1/H$. 
We can see an exponential decay of errors for every coarse mesh size $H$. For smaller $H$, the convergence is faster. This can be understood as a finite-resolution effect. For example, when $H=1/128$, there are only $H/h-1=7$ total degrees of freedom on each edge, so of course, $m=6$ basis per edge would result in a very accurate solution.

\subsection{An example with high contrast channels}
\label{sec: exp high contrast}
\begin{figure}[ht]
    \centering
    \includegraphics[width=6cm]{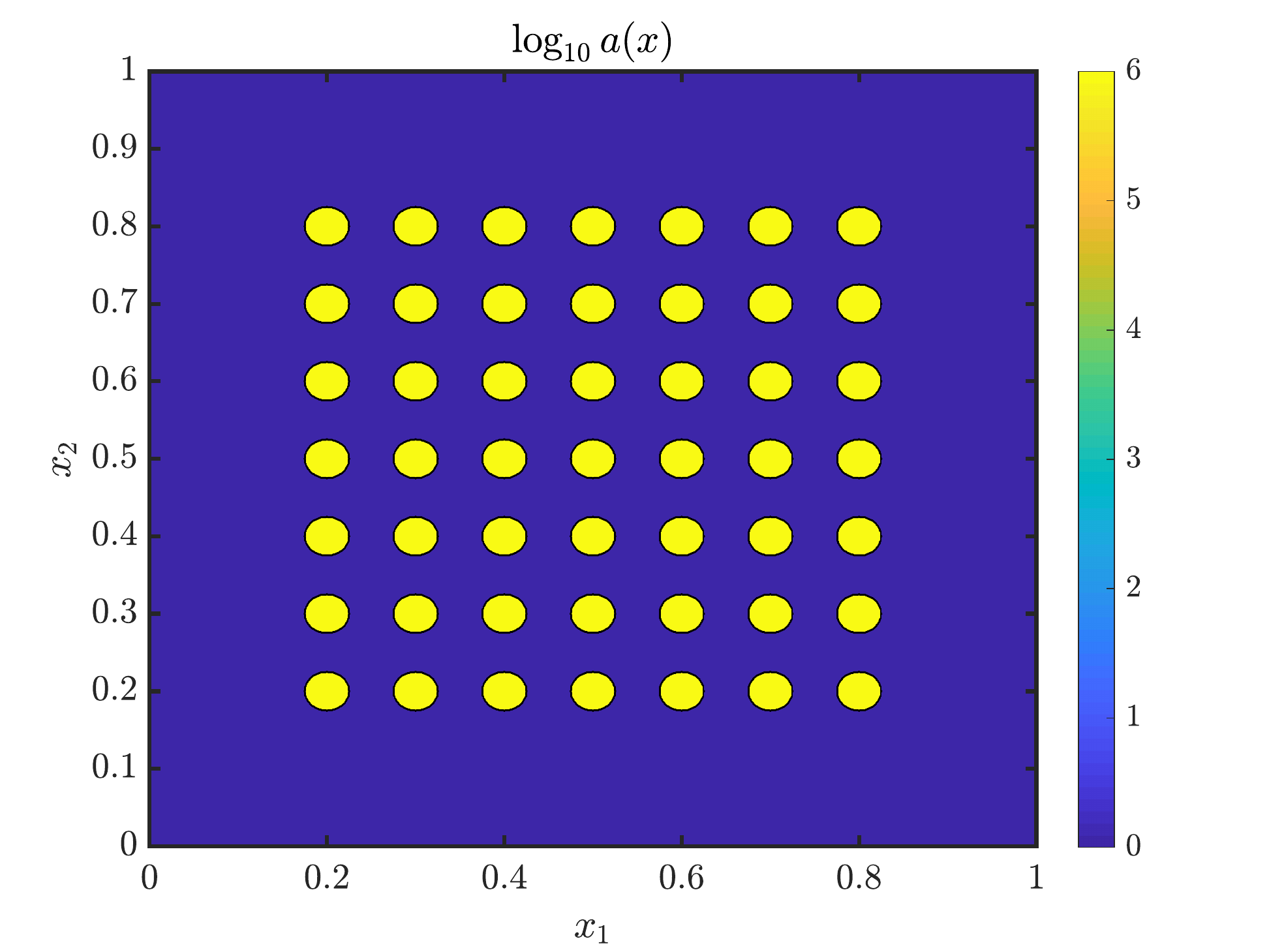}
    \includegraphics[width=6cm]{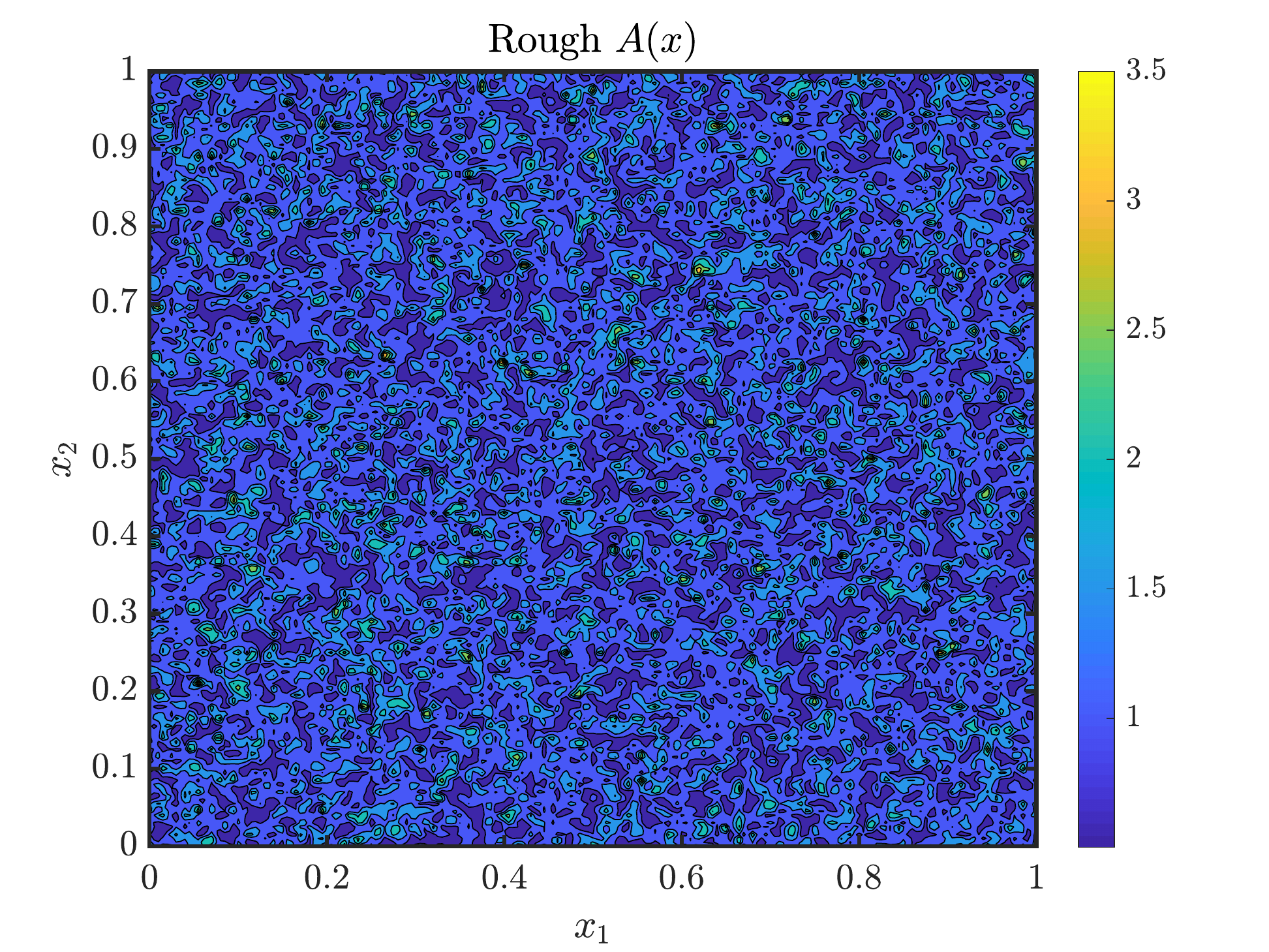}
    \caption{Left: the contour of $\log_{10} A$ for the high contrast example; right: the contour of $A$ for the rough media example.}
    \label{fig:contour_A}
\end{figure}

In the second example, we consider an elliptic problem ($V=0$) with high contrast channels. Let 
\[X:=\{(x_1,x_2) \in [0,1]^2, x_1,x_2 \in \{0.2,0.3,...,0.8\}\} \subset [0,1]^2 \, , \]
and the coefficient is defined as
\begin{equation*}
   A(x)=\left\{ 
    \begin{aligned} 1, \quad &\text{if} \ \dist(x,X)\geq 0.015\\
    M, \quad &\text{else}\, .
    \end{aligned}
    \right.
\end{equation*}
Here, $M$ is a parameter controlling the contrast. We visualize $\log_{10} A$ in the left plot of  Figure $\ref{fig:contour_A}$ for $M=10^6$.
Again, we choose homogeneous Dirichlet boundary conditions, i.e., $\Gamma_2=\emptyset$, with a non-constant right-hand side $f(x)=x_1^4-x_2^3+1$. 

\begin{figure}[ht]
    \centering
    \includegraphics[width=6cm]{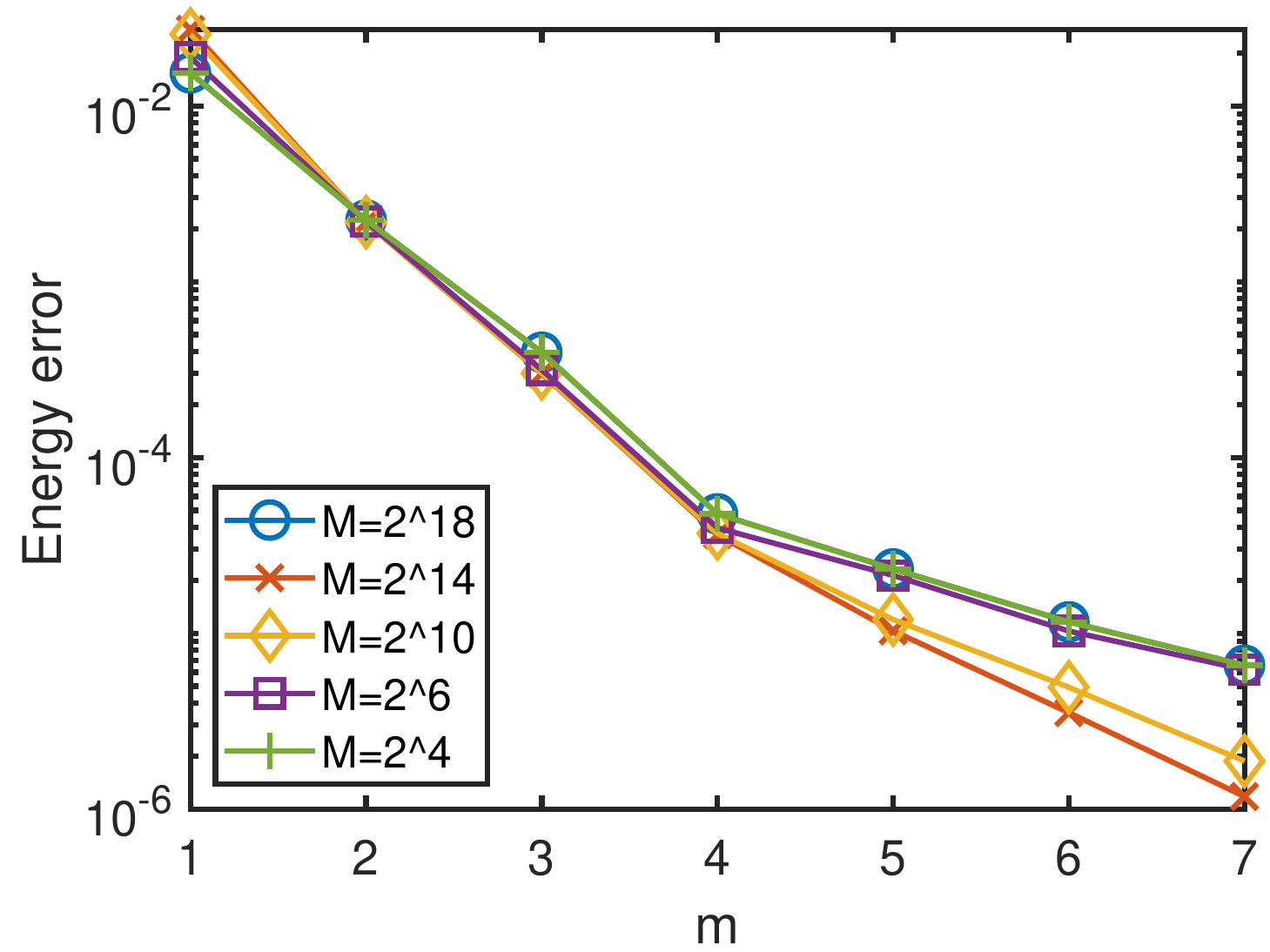}
    \includegraphics[width=6cm]{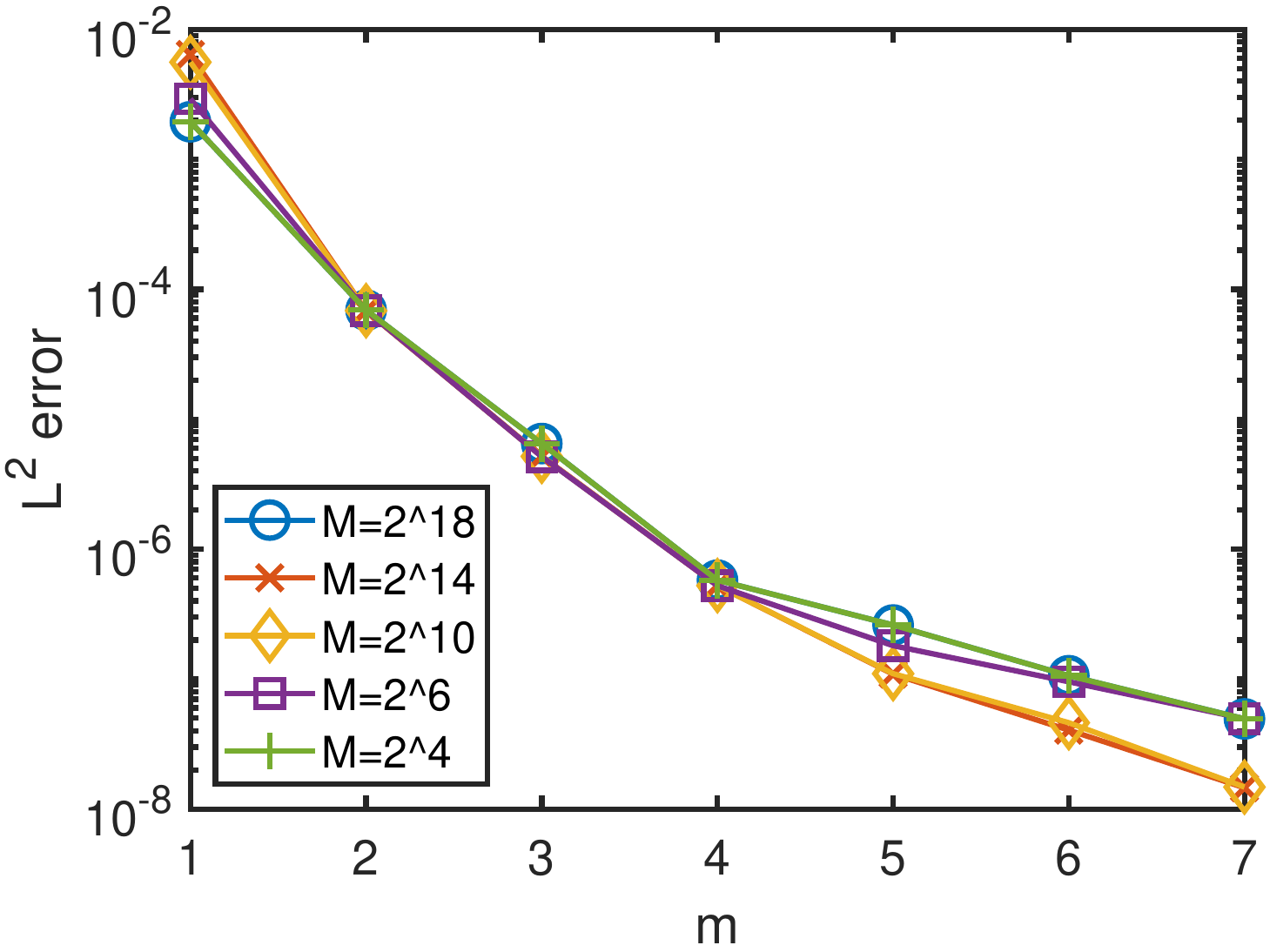}
    \caption{Numerical results for the high contrast example. Left: $e_{\cH}$ versus $m$; right: $e_{L^2}$ versus $m$.}
    \label{fig:eg2}
\end{figure}

In this example, we illustrate the convergence rate w.r.t the contrast $M$. We take different $M$ using the coarse mesh size $H=2^{-5}$ and  $m=1,2,...,7$. The numerical results are shown in Figure \ref{fig:eg2}. 
We observe a consistently exponential error decay independent of the contrast. Thus, our method demonstrates robustness with respect to the contrast $A(x)$. An intuitive explanation for this robustness could be that every step in ExpMsFEM is adaptive to $A(x)$. For example, the singular value decay of the operator $Q_{E_H}R_e$ would have some robustness regarding high contrasts in $A(x)$ because both of the norms in the domain and image of the operator is $A(x)$-weighted. We leave the theoretical analysis of deriving an $A(x)$-adapted estimates for future study. 

Also, we would like to mention that the size $h=1/1024$ of the fine mesh can actually resolve contrasts $M=2^4$ and $2^6$ only; for higher contrast, a posterior error analysis shows the reference solution on the fine mesh is not very accurate. However, we consistently observe a small error in our solution compared to the fine mesh solution, even in the regime where the fine mesh solution itself is not accurate. This implies that ExpMsFEM admits a very accurate dimension reduction of the equation on the fine mesh.

\subsection{An example of Helmholtz equation with rough field and mixed boundary}
\label{sec: exp Helmholtz}
In the last example, we consider the Helmholtz equation. This example is the same as Example 3 in \cite{chen2021exponentially}. We present it here to demonstrate that our methods are effective for complicated coefficients and mixed boundary conditions.

We impose the
homogeneous Dirichlet boundary condition on $(x_1,0), x_1 \in [0,1]$, the homogeneous Neumann boundary condition on $(x_1,1), x_1 \in [0,1]$, and the homogeneous Robin boundary condition on the other two parts of $\partial \Omega$. 
We choose $A(x)$ to be a realization of some random field; more precisely, we set
\begin{equation}
A(x)=|\xi(x)|+0.5\, ,
\end{equation}
        where the field $\xi(x)$ satisfies \[\xi(x)=a_{11}\xi_{i,j}+a_{21}\xi_{i+1,j}+a_{12}\xi_{i,j+1}+a_{22}\xi_{i+1,j+1}, \ \text{if}\ x \in [\frac{i}{2^{7}},\frac{i+1}{2^{7}})\times [\frac{j}{2^{7}},\frac{j+1}{2^{7}})\, .\]
        Here, $\{\xi_{i,j}, 0\leq i,j \leq 2^7 \}$ are i.i.d. standard Gaussian random variables. In addition, $a_{11}=(i+1-2^7x_1)(j+1-2^7x_2)$, $a_{21}=(2^7x_1-i)(j+1-2^7x_2)$, $a_{12}=(i+1-2^7x_1)(2^7x_2-j)$, $a_{22}=(2^7x_1-i)(2^7x_2-j)$ are interpolating coefficients to make $\xi(x)$ piecewise linear. A sample from this field is displayed in the right plot of Figure \ref{fig:contour_A}.

\begin{figure}[ht]
    \centering
    \includegraphics[width=6cm]{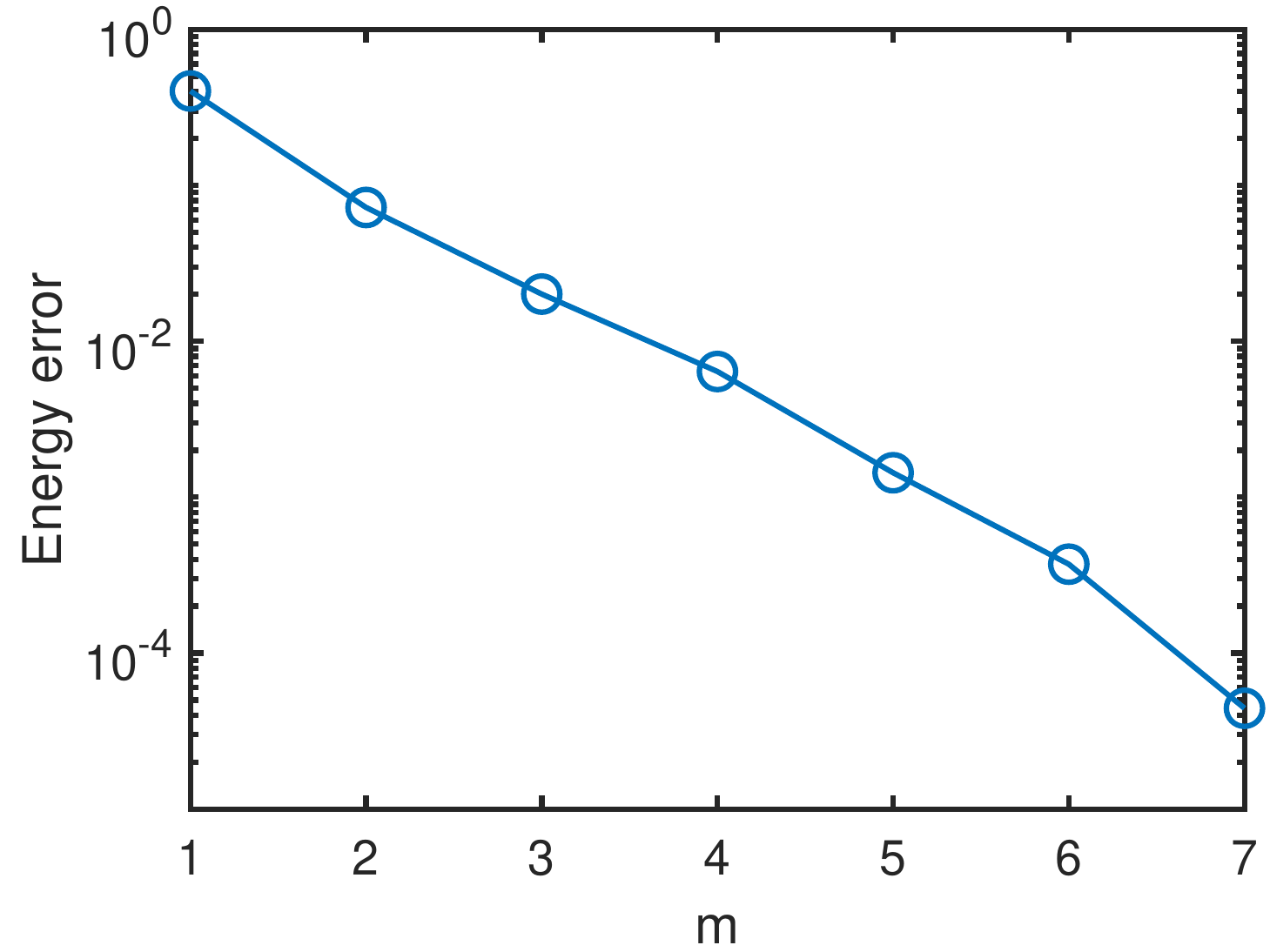}
    \includegraphics[width=6cm]{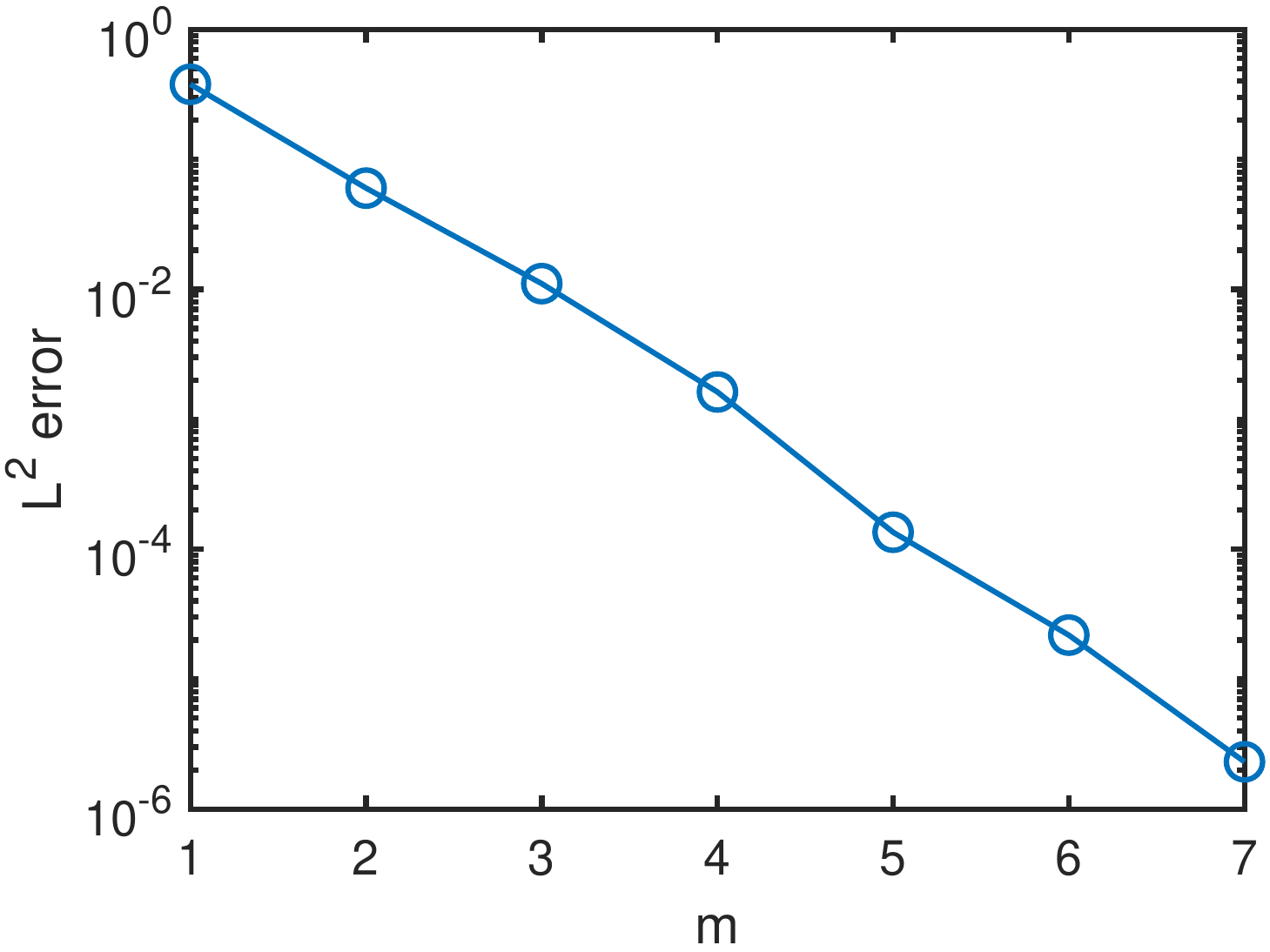}
    \caption{Numerical results for the mixed boundary and rough field example. Left: $e_{\cH}$ versus $m$; right: $e_{L^2}$ versus $m$.}
    \label{fig:eg4}
\end{figure}

Moreover, we also take $V/k^2$ and $\beta/ik$ as independent samples drawn from this random field. We choose the wavenumber $k=2^5$, the right-hand side 
{$f(x_1,x_2)=x_1^4-x_2^3+1$}, and the coarse mesh $H=2^{-5}$. Again, we take $m=1,2,...,7$ and present the numerical results in Figure \ref{fig:eg4}. 
Clearly, a nearly exponential rate of convergence is still observed for this challenging example.

\section{Discussions}
\label{sec: Discussions}
In this section, we discuss related multiscale methods in the literature; for a more specific review under the context of the elliptic and Helmholtz equations, see \cite{chen2020exponential,chen2021exponentially}. We also outline future possibilities and open questions about ExpMsFEM at the end of this section.
\subsection{Related literature}
There is a vast amount of literature on multiscale methods and numerical homogenization. 

Earlier work mainly focuses on structured $A(x)$ such as in periodic media and with scale separation; some examples include the generalized finite element methods (GFEM) \cite{babuvska1983generalized}, the multiscale finite element method (MsFEM) \cite{hou_multiscale_1997, hou1999convergence, efendiev2000convergence}, the variational multiscale methods (VMS) \cite{hughes_variational_1998}, and the heterogeneous multiscale method (HMM) \cite{abdulle2012heterogeneous}. 

Later on, people are interested in multiscale methods that can address more general rough coefficients that lie in $L^{\infty}(\Omega)$ only; see, for example, the work of optimal basis using partition of unity functions \cite{babuska2011optimal,babuvska2020multiscale, ma2021error, ma2021novel}, harmonic coordinates \cite{owhadi2007metric}, local orthogonal decomposition (LOD) \cite{malqvist_localization_2014,henning2013oversampling,kornhuber2018analysis, hauck2021super,maier2021high}, Gamblets related approaches \cite{owhadi2011localized,owhadi2014polyharmonic,owhadi2015bayesian,owhadi_multigrid_2017,hou_sparse_2017,owhadi2019operator, chen2020multiscale}, and generalizations of MsFEM \cite{hou2015optimal,chung2018constraint,li2019convergence,fu2019edge}. Different methods differ in how to find an accurate function representation. In deriving the function representation in ExpMsFEM, the solution is first decomposed into a harmonic part and a bubble part. For elliptic equations, this decomposition is the same as the orthogonal decomposition in previous work of MsFEM \cite{hou2015optimal} and approximate component mode synthesis \cite{hetmaniuk2010special,hetmaniuk2014error}.

To the best of our knowledge, among all the previous work, the optimal basis framework using partition of unity functions (and its variant) is the only one that achieves nearly exponential accuracy regarding the number of basis functions. Our ExpMsFEM \cite{chen2020exponential,chen2021exponentially} is motivated by the argument of Caccioppoli's inequality used in the optimal basis framework. ExpMsFEM is the first framework that achieves exponential accuracy without using partition of unity functions and is a direct generalization of MsFEM.

We comment in more detail on the differences and similarities between the optimal basis framework and ExpMsFEM. In the optimal basis framework, the exponentially accurate representation is obtained through the partition of unity functions rather than the edge localization and coupling in ExpMsFEM. More precisely, one can write 
\begin{equation}
\label{eqn: representation PUM}
    u = \sum_{i} \eta_i u = \sum_{i} \eta_i u^\sfh_{\omega_i} + \sum_{i} \eta_i u^\sfb_{\omega_i}\, ,
\end{equation}
where $\{\eta_i\}_i$ are partition of unity functions subordinate to an overlapped domain decomposition $\{\omega_i\}_i$ and $u_{\omega_i}^\sfh, u_{\omega_i}^\sfb$ are obtained by the harmonic-bubble splitting in $\omega_i$. The part $\eta_i u^\sfh_{\omega_i}$ can be seen as a ``restriction'' of harmonic-type functions. Thus, the argument using Caccioppoli's inequality implies that this part can be approximated by basis functions with a nearly exponential convergence rate.

Compared to \eqref{eqn: exp representation}, the representation \eqref{eqn: representation PUM} admits better geometric flexibility since by using partition of unity functions, such representation can work for problems in general dimensions. The representation \eqref{eqn: exp representation} produced by ExpMsFEM is tied to the mesh structure. When $d=2$, we have nodal and edge basis functions in the representation \eqref{eqn: exp representation}. When $d \geq 3$, we need facial basis functions and so on to represent the solution; for details see section 7 in \cite{chen2021exponentially}. In this sense, ExpMsFEM removes the partition of unity functions in the overlapped domain decomposition but pays the design cost of using a more complicated geometric structure in the non-overlapped domain decomposition. Nevertheless, the benefit of non-overlapped domain decomposition is that the basis functions are more localized since the local domain is smaller. Also, ExpMsFEM does not have the additional parameter of the partition of unity functions. Some basic numerical comparisons between ExpMsFEM and optimal basis using partition of unity functions are presented in \cite{chen2021exponentially}. We need a more in-depth comparison between the two approaches to identify their trade-offs more clearly.

\subsection{Future directions} To now, ExpMsFEM has been successfully applied to solve elliptic and Helmholtz equations. Moving forward, one can extend this idea to advection-dominated diffusion problems, time-dependent problems such as Schr\"odinger's equations, and many other linear equations. Extension to nonlinear equations appears to be nontrivial since the decomposition used in ExpMsFEM requires linearity of the equation. It could be interesting to explore the combination of ExpMsFEM and linearization to provide nonlinear homogenization of these equations.

For the current ExpMsFEM framework, we observe its robustness regarding the high contrast in the media numerically (subsection \ref{sec: exp high contrast}), but a rigorous understanding of such robustness is still lacking. Moreover, a discrete-level analysis of ExpMsFEM could be helpful for its practical use.

In essence, both ExpMsFEM and optimal basis using partition of unity functions take advantage of the low approximation complexity structures of the restriction operator on harmonic-type functions. Finding other novel low complexity structures is crucial to advance multiscale computation and model reduction. 

ExpMsFEM and optimal basis using partition of unity functions imply that nonlinear model reduction can break the Kolmogorov barrier and achieve remarkable exponential convergence. Embedding this idea to data-driven model reduction or operator learning also represents an exciting avenue for future work.

\section{Declarations}
\noindent
    {\bf Funding}. This research is in part supported by NSF Grants DMS-1912654 and DMS 2205590. We would also like to acknowledge the generous support from Mr. K. C. Choi through the Choi Family Gift Fund.

The authors have no other relevant financial or non-financial interests to disclose.
\bibliographystyle{plain}
\bibliography{ref}
\end{document}